\documentclass[12pt]{article}
\usepackage{amssymb, amsmath, amsfonts, tikz, geometry, ytableau}
\usepackage[indent,margin=1cm]{caption}
\usepackage[colorlinks,linkcolor=blue,anchorcolor=blue,citecolor=blue]{hyperref}
\usepackage{diagbox}

\numberwithin{equation}{section}
\numberwithin{figure}{section}

\hypersetup{colorlinks = true}
\newtheorem{thm}{Theorem}[section]
\newtheorem{conj}[thm]{Conjecture}

\newtheorem{lem}[thm]{Lemma}

\def\pf{\noindent{\it Proof.} }
\def\qed{\nopagebreak\hfill{\rule{4pt}{7pt}}
	\medbreak}

\allowdisplaybreaks

\def\pf{\noindent{\it Proof.} }
\def\qed{\nopagebreak\hfill{\rule{4pt}{7pt}}
	\medbreak}

\allowdisplaybreaks

\begin{document}
	\begin{center}
		{\bf \Large Interlacing property of a family of generating polynomials over Dyck paths}
	\end{center}
\begin{center}
	{\bf Bo Wang$^1$ and Candice X.T. Zhang$^2$\\[6pt]}
	
	{\it Center for Combinatorics, LPMC\\
		Nankai University, Tianjin 300071, P. R. China\\[8pt]
		
		Email: $^{1}${\tt bowang@nankai.edu.cn}, $^{2}${\tt zhang\_xutong@mail.nankai.edu.cn}}
\end{center}

\noindent\textbf{Abstract.}
In the study of a tantalizing symmetry on Catalan objects, B\'ona et al. introduced a family of polynomials $\{W_{n,k}(x)\}_{n\geq k\geq 0}$ defined by
\begin{align*}
W_{n,k}(x)=\sum_{m=0}^{k}w_{n,k,m}x^{m},
\end{align*}
where $w_{n,k,m}$ counts the number of Dyck paths of semilength $n$ with $k$ occurrences of $UD$ and $m$ occurrences of $UUD$.
They proposed two conjectures on the interlacing property of these polynomials, one of which states that $\{W_{n,k}(x)\}_{n\geq k}$ is a Sturm sequence for any fixed $k\geq 1$, and the other states that $\{W_{n,k}(x)\}_{1\leq k\leq n}$ is a Sturm-unimodal sequence for any fixed $n\geq 1$. In this paper, we obtain certain recurrence relations for $W_{n,k}(x)$, and further confirm their conjectures.

\noindent \textbf{AMS Classification 2020:} 05A15, 26C10


\section{Introduction}
A {\it Dyck path} of semilength $n$ in $\mathbb{Z}^2$ is a lattice path starting at the origin $(0,0)$, ending at $(2n,0)$, and never goes below the $x$-axis, whose permitted step types are up step $U=(1,1)$ and down step $D=(1,-1)$. It is well known that the set of Dyck paths of semilength $n$ is counted by the Catalan number $C_{n}=\frac{1}{n+1}\binom{2n}{n}$, which is the sequence A000108 in the On-line Encyclopedia of Integer Sequences of Sloane \cite{Sloane-2018}.

Numerous studies have been focused on the refinement of Catalan numbers by considering certain statistics over Dyck paths. It is easy to see that a Dyck path determines a word of $\{U,D\}$ as one records the steps along the path from left to right. One important class of statistics are defined with various factors appearing in the word representation of Dyck paths. It seems that the most natural factor is a $UD$-factor, which means that an up step is immediately followed by a down step in the Dyck path. The number of Dyck paths of semilength $n$ with exactly $k$ occurrences of $UD$-factors is given by the Narayana number $N(n,k)=\frac{1}{n}\binom{n}{k-1}\binom{n}{k}$; see Sulanke \cite{Sulanke-1999}.
The enumeration of Dyck paths of semilength $n$ with $k$ occurrences of $UUD$-factors has been first stuided by Sapounakis, Tasoulas, and Tsikouras \cite{Sapounakis-Tasoulas-Tsikouras-2006}.
Lin and Kim \cite{Lin-Kim-2021} introduced the segment statistic, which is actually the $UUD$-factor, to study various classical statistics on restricted inversion sequences. In their paper, Lin and Kim also proved that this statistic, when applied to Dyck paths, is equidistributed with the descent statistic over the group of $(3, 2, 1)$-avoiding permutations.
Wang \cite{Wang-2011} developed a useful technique for computing relevant generating functions for Dyck paths with different factors. For more information on the enumeration of Dyck paths with respect to various factors, see \cite{Denise-1995,Deutsch-1999,Merlini-Sprugnoli-Verri-2002,Sun-2004,Woan-2004,Mansour-2006,Czabarka-Florez-2015} and references therein.

This paper is much motivated by a recent work \cite{Bona-Dimitrov-2022} due to B\'{o}na et al., who first considered the joint distribution of $UD$-factors and $UUD$-factors over Dyck paths.
Let $w_{n,k,m}$ be the number of Dyck paths of semilength $n$ with $k$ $UD$-factors and $m$ $UUD$-factors. For these numbers $w_{n,k,m}$, B\'{o}na et al. \cite{Bona-Dimitrov-2022} proved the following tantalizing symmetric property:
\[ w_{2k+1,k,m}=w_{2k+1,k,k+1-m}, \quad \mbox{where } 1\le m\le k.\]
To obtain this result,
they derived the following explicit formula for the numbers $w_{n,k,m}$ by using generating function techniques.

\begin{thm}
	[\cite{Bona-Dimitrov-2022}, Theorem 1.2] For all $n$, $k$ and $m$, we have
	\begin{equation}\label{eq-wnkm}
		w_{n,k,m}=
		\begin{cases}
			\frac{1}{k}\binom{n}{k-1}\binom{n-k-1}{m-1}\binom{k}{m}, & \mbox{if } 0<m\leq k, \mbox{and } k+m\leq n,\\
			1, & \mbox{if } m=0\ \mbox{and } n=k,\\
			0, & \mbox{otherwise.}\\
		\end{cases}
	\end{equation}
\end{thm}
With the above formula, B\'{o}na et al. also noted that the numbers $w_{n,k,m}$
are closely related to the classical Narayana Numbers, as well as to Callan's generalization of Narayana Numbers \cite{Callan-2017}.

Let $W_{n,k}(x)$ be the generating polynomial of $w_{n,k,m}$ as given by
\begin{equation}\label{eq-wnkx}
	W_{n,k}(x)=\sum_{m=0}^{k}w_{n,k,m}x^{m}.
\end{equation}
By \eqref{eq-wnkm}, it is clear that $\deg W_{n,k}(x)=\min \{k,n-k \}$ if $n> k$ and $\deg W_{n,k}(x)=0$ otherwise.
It turns out that these polynomials enjoy very interesting properties. B\'{o}na et al. \cite{Bona-Dimitrov-2022} obtained the following result. 

\begin{thm}[\cite{Bona-Dimitrov-2022}, {Proposition 6.1, Theorem 6.3}]
	\label{same root}
For any $n,k\geq 0$ the polynomial $W_{n,k}(x)$ has only real zeros. Moreover, for all $1\leq k\leq n-1$, the polynomials $W_{n,k}(x)$ and $W_{n,n-k}(x)$ have the same zeros.
\end{thm}

The real-rootedness of $W_{n,k}(x)$ was proved by B\'{o}na et al. \cite{Bona-Dimitrov-2022} based on Malo’s result regarding the roots of
the Hadamard product of two real-rooted polynomials.
It is worth mentioning that many combinatorial polynomials have only real zeros. For excellent surveys on this topic, we refer the readers to Stanley  \cite{Stanley-log-1989}, Brenti \cite{Brenti-1994}, and Br\"and\'en \cite{Branden-2015}. One of the useful methods to prove the real-rootedness of a polynomial is to consider the interlacing property involving its zeros. 

B\'{o}na et al. \cite{Bona-Dimitrov-2022} further studied the interlacing property of $W_{n,k}(x)$ by fixing $n$ or $k$, and proposed two interesting conjectures.
Before stating these conjectures, let us recall some related definitions following Liu and Wang \cite{Liu-Wang-2006}. Given two real-rooted polynomials $F(x)$ and $G(x)$ with nonnegative real coefficients, let $\{\alpha_{r}\}$ and $\{\beta_{s}\}$ be their zeros in weakly decreasing order, respectively. We say that $G(x)$ {\it interlaces} $F(x)$, denoted by $G(x)\preccurlyeq F(x)$, if
$\deg F(x)=\deg G(x)=n$ and 
\begin{equation*}\label{ineq-interlace-1}
\beta_{n}\leq \alpha_{n}\leq \beta_{n-1}\leq \alpha_{n-1}\leq \cdots\leq \beta_{1}\leq \alpha_{1},
\end{equation*}
or $\deg f(x)=\deg g(x)+1=n$ and
\begin{equation*}\label{ineq-interlace-2}
	\alpha_{n}\leq \beta_{n-1}\leq \alpha_{n-1}\leq \cdots\leq \beta_{1}\leq \alpha_{1}.
\end{equation*}
For convenience, let $a\preccurlyeq bx+c$ for any real numbers $a,b,c$ and $F(x)\preccurlyeq 0$,  $0\preccurlyeq F(x)$ for any real-rooted polynomial $F(x)$.
Given a sequence $\{F_{i}(x)\}_{i\geq 0}$ of real-rooted polynomials, we say that it is a {\it generalized Sturm sequence} if $F_{i}(x)\preccurlyeq F_{i+1}(x)$ for all $i\geq 0$. We would like to point out that a generalized Sturm sequence is called a Sturm sequence by B\'{o}na et al. They also introduced the notion of Sturm-unimodal sequences. With our notation here, a finite sequence $\{F_{i}(x)\}_{1\leq i\leq n}$ of real-rooted polynomials is said to be {\it Sturm-unimodal}, provided that there exists $1\leq j\leq n$ such that
\[F_{1}(x)\preccurlyeq \cdots\preccurlyeq F_{j-1}(x)\preccurlyeq F_{j}(x)\succcurlyeq F_{j+1}(x)\succcurlyeq \cdots\succcurlyeq F_{n}(x).\]
Now the two conjectures of B\'{o}na et al. \cite{Bona-Dimitrov-2022} can be stated as follows. 

\begin{conj}[\cite{Bona-Dimitrov-2022}, Conjecture 6.4]\label{conje1}
 For any fixed $k\geq 1$, 
 the polynomial sequence $\{W_{n,k}(x)\}_{n\geq k}$
 is a generalized Sturm sequence.
\end{conj}

\begin{conj}[\cite{Bona-Dimitrov-2022}, Conjecture 6.5]\label{conje2}
 For any fixed $n\geq 1$, the polynomial sequence $\{W_{n,k}(x)\}_{1\leq k\leq n}$ is Sturm-unimdoal.
\end{conj}

In this paper we shall prove these two conjectures.

\section{The main results}

The aim of this section is to prove Conjecture \ref{conje1} and Conjecture \ref{conje2}. In the process of proving these two conjectures, 
we need the following result which provides a sufficient condition for a polynomial sequence with three-term recurrence to be a generalized Sturm sequence. Note that it is a special case of Liu and Wang's criterion \cite{Liu-Wang-2006}.

\begin{thm}	[\cite{Liu-Wang-2006}, Corollary 2.4]\label{thm-lw}
	Let $\{F_{i}(x)\}_{i\geq 0}$ be a sequence of polynomial with nonnegative coefficients satisfying the following conditions:
	\begin{itemize}
		\item[1.] $F_{0}(x)$ and $F_{1}(x)$ are real-rooted polynomial with $F_{0}(x)\preccurlyeq F_{1}(x)$.
		
		\item[2.] $\deg F_{i+1}(x)=\deg F_{i}(x)$ or $\deg F_{i}(x)+1$ for any $i\geq 0$.
		
		\item[3.] There exist polynomials $A_{j}(x)$ and $B_{j}(x)$ with real coefficients such that
		\begin{align}\label{eq-formal-rec}
		F_{j+2}(x)=A_{j}(x)F_{j+1}(x)+B_{j}(x)F_{j}(x).
		\end{align}
	\end{itemize}
	If for all $x\leq 0$, we have  $B_{j}(x)\leq 0$, then $\{F_{i}(x)\}_{i\geq 0}$ is a generalized Sturm sequence.
\end{thm}

In order to use the above theorem to prove Conjecture \ref{conje1}, we need to 
establish some recurrence relation satisfied by the polynomials $W_{n,k}(x)$ when fixing $k$. 
Let us first give a recurrence relation of the coefficients $w_{n,k,m}$ for each fixed $k\geq 1$.

\begin{lem}\label{thm-conje1-wnkm-recu}
	Let $w_{n,k,m}$ be as given by \eqref{eq-wnkm}. Then, for any $n\geq k-1\geq 0$, we have 
\begin{equation}\label{eq-thm-conje1-wnkm-recu}
\begin{split}
w_{n+2,k,m}=&\frac{2(n+2)(n-k+1)}{(n-k+2)(n-k+3)}w_{n+1,k,m}-\frac{(n+2)(n-2k+1)}{(n-k+2)(n-k+3)}w_{n+1,k,m-1}\\
&+\frac{(n+1)(n+2)(n-k)}{(n-k+2)^2(n-k+3)}(w_{n,k,m-1}-w_{n,k,m}).
\end{split}
\end{equation}
\end{lem}

\pf We may assume that $0\le m\leq k$ and $n\geq m+k-2$, since there is nothing to prove for $m<0$, or $m>k$, or $n< m+k-2$.
Moreover, it is routine to verify the validity of \eqref{eq-thm-conje1-wnkm-recu} for $m=0$ since both sides vanish under the condition $n\ge k-1$.

If $m=1$, and hence $n\ge k-1$, then we divide the proof of \eqref{eq-thm-conje1-wnkm-recu} into the following three cases.

{\bf Case A1:} $n=k-1$. We find that 
\[w_{n+1,k,m}=w_{n,k,m-1}=w_{n,k,m}=0\text{ and }w_{n+1,k,m-1}=1\]
in view of \eqref{eq-wnkm}. Thus, \eqref{eq-thm-conje1-wnkm-recu} holds since its right-hand side simplifies to $\binom{k+1}{2}$, which is indeed equal to $w_{n+2,k,m}$. 

{\bf Case A2:} $n=k$. In this case, the third term on the right-hand side of \eqref{eq-thm-conje1-wnkm-recu} naturally vanishes.  
Note that for $m=1$ we have $w_{n+1,k,m-1}=0$ by \eqref{eq-wnkm}. 
Thus, it suffices to show that
\[w_{k+2,k,1}=\frac{2(k+2)}{3}w_{k+1,k,1}, \]
which can be easily verified by \eqref{eq-wnkm}.

{\bf Case A3:}  $n>k$. Keep in mind that $m=1$ throughout this case. By \eqref{eq-wnkm} we have $w_{n+1,k,m-1}=w_{n,k,m-1}=0$. Thus, the right-hand side of \eqref{eq-thm-conje1-wnkm-recu} turns out to be
\[\frac{2(n+2)(n-k+1)}{(n-k+2)(n-k+3)}\binom{n+1}{k-1}-\frac{(n+1)(n+2)(n-k)}{(n-k+2)^2(n-k+3)}\binom{n}{k-1}=\binom{n+2}{k-1},\] 
which is equal to $w_{n+2,k,m}$, as desired.

From now on we may assume that $2\leq m\leq k$.
We further divide the proof of \eqref{eq-thm-conje1-wnkm-recu} into the following three cases.

{\bf Case B1:} $n=m+k-2$. In this case, we have $w_{n+1,k,m}=w_{n,k,m-1}=w_{n,k,m}=0$.
Now it is sufficient to show that
\[w_{n+2,k,m}=-\frac{(n+2)(n-2k+1)}{(n-k+2)(n-k+3)}w_{n+1,k,m-1},\]
which can be verified by \eqref{eq-wnkm}.

{\bf Case B2:} $n=m+k-1$. For this case we have $w_{n,k,m}=0$ by \eqref{eq-wnkm}.
Subtituting \eqref{eq-wnkm} and the condition $m=n-k+1$ into the right-hand side of \eqref{eq-thm-conje1-wnkm-recu} and then simplifying, we obtain that
\begin{align*}
&\frac{n-k+1}{k(n-k+2)}\binom{n+2}{k-1}\binom{k}{n-k+1}\binom{n-k}{n-k-1}+\frac{2(n-k+1)}{k(n-k+2)}\binom{n+2}{k-1}\binom{k}{n-k+1}\\[5pt]
&=\frac{1}{k}(n-k+1)\binom{n+2}{k-1}\binom{k}{n-k+1},
\end{align*}
which is equal to $w_{n+2,k,m}$ according to \eqref{eq-wnkm}.

{\bf Case B3:} $n\geq m+k$.  
By substituting \eqref{eq-wnkm} into the right-hand side of \eqref{eq-thm-conje1-wnkm-recu} and then simplifying, we get
\begin{equation*}
	\begin{split}	
		&\frac{n-k+1+m}{k(n-k+2)}\binom{n+2}{k-1}\binom{k}{m}\binom{n-k}{m-1}+\frac{m}{k(n-k+2)}\binom{n+2}{k-1}\binom{k}{m}\binom{n-k}{m-2}\\[5pt]
		&=\frac{1}{k}\binom{n+2}{k-1}\binom{k}{m}\left(\frac{n-k+m+1}{n-k+2}\binom{n-k}{m-1}+\frac{m}{n-k+2}\binom{n-k}{m-2}\right)\\[5pt]
		&=\frac{1}{k}\binom{n+2}{k-1}\binom{n-k+1}{m-1}\binom{k}{m},
	\end{split}
\end{equation*}
which is equal to $w_{n+2,k,m}$ according to \eqref{eq-wnkm}, as desired. 

Taking into account all the above cases, we complete the proof.  
\qed

The recurrence relation \eqref{eq-thm-conje1-wnkm-recu} satisfied by the coefficients $w_{n,k,m}$ is equivalent to the following recurrence relation for polynomials $W_{n,k}(x)$, which plays a key role in our proof of Conjecture \ref{conje1}.

\begin{thm}\label{thm-conje1-recu}
Fixing $k\geq 1$, for any  $n\ge k-1$ we have
\begin{equation}\label{eq-thm-conje1-recu}
\begin{split}
W_{n+2,k}(x)=&\frac{(n+2)\left(2(n-k+1)-(n-2k+1)x\right)}{(n-k+2)(n-k+3)}W_{n+1,k}(x)\\[5pt]
&+\frac{(n+1)(n+2)(n-k)(x-1)}{(n-k+2)^2(n-k+3)}W_{n,k}(x).
\end{split}
\end{equation}
\end{thm}

\pf Due to the fact that $\deg W_{n,k}(x)=\min\{k,n-k\}$, it suffices to compare the coefficients of $x^m$ for $0\leq m\leq k$.  
For the right-hand side of \eqref{eq-thm-conje1-recu}, the coefficient of $x^m$ is 
\begin{equation*}\label{eq-proof-conje1-2}
\begin{split}
&\frac{2(n+2)(n-k+1)}{(n-k+2)(n-k+3)}w_{n+1,k,m}-\frac{(n+2)(n-2k+1)}{(n-k+2)(n-k+3)}w_{n+1,k,m-1}\\
&\qquad\qquad+\frac{(n+1)(n+2)(n-k)}{(n-k+2)^2(n-k+3)}(w_{n,k,m-1}-w_{n,k,m}),
\end{split}
\end{equation*}
which is equal to $w_{n+2,k,m}$ by \eqref{eq-thm-conje1-wnkm-recu}. 
This is just the coefficent  of $x^m$ on the left-hand side of \eqref{eq-thm-conje1-recu}. The proof is complete. 
\qed

We proceed to prove Conjecture \ref{conje1}. Using Theorem \ref{thm-conje1-recu} and Theorem \ref{thm-lw}, we obtain the first main result of this section.

\begin{thm}\label{thm-conj1}
For any fixed $k\geq 1$, the polynomial sequence $\{W_{n,k}(x)\}_{n\geq k}$ is a generalized Sturm sequence.
\end{thm}

\pf
Taking $F_{i}(x)$ in Theorem \ref{thm-lw} to be the polynomial 
$W_{k+i,k}(x)$ for each $i\geq 0$, it is clear that each $F_{i}(x)$ is a polynomial with nonnegative coefficients. By \eqref{eq-wnkm} and \eqref{eq-wnkx}, we have $\deg F_i(x)=i$ if $i\le k$ and $\deg F_i(x)=k$ otherwise.
Note that
$$F_{0}(x)=W_{k,k}(x)=1,\ \ \ \ F_{1}(x)=W_{k+1,k}(x)=\dbinom{k+1}{k-1}x,$$
and hence $F_{0}(x)\preccurlyeq F_{1}(x)$. 
Now the recurrence relation (\ref{eq-thm-conje1-recu}) can be restated as in the form \eqref{eq-formal-rec} with
\begin{equation*}
	A_{j}(x)=\frac{(k+j+2)}{(j+2)(j+3)}\left(2(j+1)-(j-k+1)x\right),
\end{equation*} 
and
\begin{equation*}
B_{j}(x)=\frac{(k+j+1)(k+j+2)j}{(j+2)^2(j+3)}(x-1).
\end{equation*}
Clearly, for any $j\geq 0$ and $x\leq 0$, we have $B_{j}(x)\leq 0$. From Theorem \ref{thm-lw} it follows that the sequence $\{F_{i}(x)\}_{i\geq 0}$, and hence $\{W_{n,k}(x)\}_{n\geq k}$, is a  generalized Sturm sequence.
\qed

In order to prove Conjecture \ref{conje2}, we present the following recurrence relation of the coefficients $w_{n,k,m}$ for each fixed $n\geq 1$.

\begin{lem}\label{thm-conje2-wnkm-recu}
Let $w_{n,k,m}$ be as given by \eqref{eq-wnkm}. Then, for any $1\le k\le \lfloor\frac{n+1}{2} \rfloor-2$, we have
\begin{equation}\label{eq-thm-conje2-wnkm-recu}
w_{n,k+2,m}=a(n,k)w_{n,k+1,m-1}+b(n,k)w_{n,k+1,m}-c(n,k)w_{n,k,m},
\end{equation}
where
\begin{equation*}
\left\{
\begin{aligned}
a(n,k)&=\frac{(n-k)(n-2k-2)(n-2k-3)}{(k+1)(k+2)(n-k-2)},\\[5pt]
b(n,k)&=\frac{2(n-k)(n-k-1)(n-2k-2)}{(k+2)(n-k-2)(n-2k-1)},\\[5pt]
c(n,k)&=\frac{(n-k+1)(n-k)^{2}(n-2k-3)}{(k+1)(k+2)(n-k-2)(n-2k-1)}.
\end{aligned}
\right.
\end{equation*}
\end{lem}
\pf
By similar arguments as in the proof of Lemma \ref{thm-conje1-wnkm-recu}, we may assume that $1\le m\leq k$, and hence $2\leq m+k\le n-2$ by the condition $1\le k\le \lfloor \frac{n+1}{2}\rfloor-2$.

If $m=1$, then we have $w_{n,k+1,m-1}=0$ by \eqref{eq-wnkm}, and the right-hand side of \eqref{eq-thm-conje2-wnkm-recu} turns out to be
\[b(n,k)\binom{n}{k}-c(n,k)\binom{n}{k-1}=\binom{n}{k+1}=w_{n,k+2,m},\]
as desired.

The proof for the case of $2\le m\le k$ is very similar to that of Case B3 in the proof of Lemma \ref{thm-conje1-wnkm-recu}, and is omitted here. 
\qed

Fixing an integer $n\geq 1$, the above lemma immediately leads to a recurrence relation for the polynomial sequence $\{W_{n,k}(x)\}_{1\leq k\leq \lfloor\frac{n+1}{2} \rfloor}$. 

\begin{thm}\label{thm-conje2-recu}
Let $a(n,k)$, $b(n,k)$ and $c(n,k)$ be as given in Lemma \ref{thm-conje2-wnkm-recu}. Then, for any $1\le k\le \lfloor\frac{n+1}{2} \rfloor-2$, we have
\begin{equation}\label{eq-thm-conje2-recu}
W_{n,k+2}(x)=\left(a(n,k)x+b(n,k)\right)W_{n,k+1}(x)-c(n,k)W_{n,k}(x).
\end{equation}
\end{thm}
\pf The desired recurrence immediately follows by comparing the coefficients of $x^{m}$ on both sides of (\ref{eq-thm-conje2-recu}) for each $m$ and then using \eqref{eq-thm-conje2-wnkm-recu}. 
\qed

Based on the recurrence relation \eqref{eq-thm-conje2-recu} and Theorem \ref{thm-lw}, we obtain the following result, which provides an affirmative answer to Conjecture \ref{conje2}.

\begin{thm}\label{thm-conj2}
For any fixed $n\geq 1$, the polynomial sequence $\{W_{n,k}(x)\}_{1\leq k\leq n}$ is Sturm-unimodal.
\end{thm}
\pf
By \eqref{eq-wnkm} and \eqref{eq-wnkx}, one can check that 
\begin{equation}\label{eq-sturm-unimodal}
\begin{split}
W_{n,n}(x)&=1, \mbox{ for } n\geq 1,\\
W_{n,1}(x)&=x, \mbox{ for } n\geq 2,\\ W_{n,2}(x)&=\frac{n}{2}\left(2x+(n-3)x^2\right), \mbox{ for } n\geq 3.
\end{split}
\end{equation}
Thus we always have $W_{n,n-1}(x)\succcurlyeq W_{n,n}(x)$ by convention.
On the other hand, Theorem \ref{same root} implies that 
$W_{n,k}(x)\preccurlyeq W_{n,k+1}(x)$ if and only if $W_{n,n-k-1}(x)\succcurlyeq W_{n,n-k}(x)$ for any $1\le k\le \lfloor\frac{n+1}{2} \rfloor$. Therefore, it suffices to show that the polynomial sequence  $\{W_{n,k}(x)\}_{1\leq k\leq \lfloor\frac{n+1}{2}\rfloor}$  is a generalized Sturm sequence. For  $1\le n\le 5$ this can be directly verified by \eqref{eq-sturm-unimodal}.  

Now suppose that $n\geq 6$ for the remainder of the proof. 
For $0\le i\le \lfloor\frac{n+1}{2}\rfloor-1$, take each $F_{i}(x)$ in Theorem \ref{thm-lw} to be the polynomial $W_{n,i+1}(x)$. Clearly, $F_{i}(x)$ is a polynomial with nonnegative coefficients. By \eqref{eq-wnkm} and \eqref{eq-wnkx}, we have
\begin{align}
\deg F_{i}(x)=\left\{
\begin{array}{ll}
i+1, & \mbox{ if } 0\le i\le \lfloor\frac{n+1}{2}\rfloor-2,\\
i+1, & \mbox{ if $n$ is even and } i=\lfloor\frac{n+1}{2}\rfloor-1,\\
i, & \mbox{ if $n$ is odd and } i=\lfloor\frac{n+1}{2}\rfloor-1.
\end{array}
\right.
\end{align} 
It is also clear that $F_{0}(x)\preccurlyeq F_{1}(x)$. 
Now the recurrence relation (\ref{eq-thm-conje2-recu}) can be restated as in the form \eqref{eq-formal-rec} with
\begin{equation*}
	A_{j}(x)=\frac{(n-j-1)(n-2j-4)(n-2j-5)}{(j+2)(j+3)(n-j-3)}x+\frac{2(n-j-1)(n-j-2)(n-2j-4)}{(i+3)(n-j-3)(n-2j-3)},
\end{equation*}
and
\begin{equation*}
	B_{j}(x)=-\frac{(n-j)(n-j-1)^{2}(n-2j-5)}{(j+2)(j+3)(n-j-3)(n-2j-3)}\leq 0
\end{equation*}
for any $0\le j\le \lfloor\frac{n+1}{2}\rfloor-3$.
We conclude from Theorem \ref{thm-lw} that the polynomial sequence $\{F_{i}(x)\}_{0\leq i\leq \lfloor\frac{n+1}{2}\rfloor-1}$, and hence $\{W_{n,k}(x)\}_{1\leq k\leq \lfloor\frac{n+1}{2}\rfloor}$, is a generalized Sturm sequence. 
This completes the proof.
\qed

\vskip 0.5cm
\noindent \textbf{Acknowledgments.}
This work was supported by the Fundamental Research Funds for the Central Universities and the National Science Foundation of China (Grant No. 11971249).

\end{document}